\documentclass[twoside,11pt]{article}
\title{} \author{} \date{}
\usepackage{latexsym,amssymb,times}
\input amssym.def
\newtheorem{te}{Theorem}[section]
\newtheorem{prop}[te]{Proposition}

\newtheorem{fac}[te]{Fact}
\newtheorem{lem}[te]{Lemma}

\newtheorem{rem}[te]{Remark}
\newtheorem{ex}[te]{Example}

\def\dok{\noindent{\bf Proof. }}
\def\kdok{\hfill $\Box$ \par \vspace*{2mm} }
\def\a{\alpha}
\def\b{\beta}
\def\g{\gamma}
\def\d{\delta}
\def\f{\varphi}
\def\p{\psi}
\def\o{\omega}
\def\k{\kappa}
\def\l{\lambda}
\def\r{\rho}
\def\s{\sigma}

\def\c{{\mathfrak c}}

\def\P{{\mathbb P}}
\def\Q{{\mathbb Q}}

\def\N{{\mathbb N}}
\def\X{{\mathbb X}}

\def\Z{{\mathbb Z}}
\def\D{{\mathbb D}}

\def\BR{{\mathbb R}}

\def\BT{{\mathbb T}}
\def\BL{{\mathbb L}}

\def\CO{{\mathcal O}}

\def\I{{\mathcal I}}

\def\CN{{\mathcal N}}
\def\F{{\mathcal F}}

\def\down{\!\downarrow}
\def\up{\!\uparrow}
\def\la{\langle}
\def\ra{\rangle}

\def\dom{\mathop{\mathrm{dom}}\nolimits}
\def\ran{\mathop{\mathrm{ran}}\nolimits}
\def\id{\mathop{\mathrm{id}}\nolimits}

\def\Lev{\mathop{\mathrm{Lev}}\nolimits}

\def\Iso{\mathop{\rm Iso}\nolimits}
\def\Aut{\mathop{\rm Aut}\nolimits}

\def\Lim{\mathop{\rm {Lim}}\nolimits}

\def\ar{\mathop{\rm ar}\nolimits}

\def\Mono{\mathop{\rm Mono}\nolimits}

\def\Sur{\mathop{\rm Sur}\nolimits}

\def\Cond{\mathop{\rm Cond}\nolimits}
\def\Sym{\mathop{\rm Sym}\nolimits}

\def\otp{\mathop{\rm otp}\nolimits}

\def\he{\mathop{\rm ht}\nolimits}
\def\Ord{\mathop{\rm Ord}\nolimits}
\def\Min{\mathop{\rm Min}\nolimits}
\def\Max{\mathop{\rm Max}\nolimits}
\def\is{\mathop{\rm Is}\nolimits}

\def\Pc{\mathop{\rm Pc}\nolimits}
\def\bcd{\mathop{\dot{\bigcup}}\nolimits}
\def\du{\mathop{\mathrel{\dot{\cup}}}\nolimits}
\def\Prt{\mathop{\rm Prt}\nolimits}
\def\Tree{\mathop{\rm Tree}\nolimits}
\begin{document}
\thispagestyle{plain}
\vspace*{-16mm}
\begin{center}
           {\large \bf {\uppercase{Reversible and irreversible trees}}}
\end{center}
\begin{center}
{\bf Milo\v s S.\ Kurili\'c}\footnote{Department of Mathematics and Informatics, Faculty of Sciences, University of Novi Sad,
                                      Trg Dositeja Obradovi\'ca 4, 21000 Novi Sad, Serbia,
                                      e-mail: milos@dmi.uns.ac.rs}
\end{center}
\begin{abstract}
\noindent
A structure $\X$ is reversible iff every bijective endomorphism of $\X$ is an automorphism.
A combinatorial characterization of that property in the class of trees is obtained: a tree $\BT$ is non-reversible
iff it contains a {\it critical node} or an {\it archetypical subtree}
(suborders of $\BT$ with some additional properties).
Consequently, a tree with finite nodes $\BT$ is reversible iff it does not contain archetypical subtrees.
In particular, $\BT$ is reversible if for each ordinal $\a \in [\o ,\he (\BT ))$
the sequence $\la|N|: \CN (\BT ) \ni N\subset L_\a\ra $ is  finite-to-one or almost constant.
So, regular $n$-ary trees are reversible,
reversible Aronszajn trees exist
and, if there are Suslin or Kurepa trees, there are reversible ones.

Using a characterization of reversibility via back and forth systems 
a wide class of non-reversible trees is detected:
{\it bad trees} (having all branches of height $\he (\BT)=|T|=|L_0|$, where $|T|$ is a regular cardinal).
Consequently, a countable tree $\BT$ of height $\o$ and without leafs is reversible iff its nodes are finite
iff the set of its branches is a compact subset of the Baire space.

The results concerning disconnected trees are relevant in this context. 
So, the tree $\bcd _\mu {}^{<\a }\l$ is reversible iff $\min\{\a ,\l\mu\} <\o$. 
If $\lambda \geq \o$, then $\bcd _{i\in I}{}^{<n_i}\l$ is a reversible tree 
iff the sequence $\la n_i :i\in I\ra $ is  finite-to-one or almost constant at some $k\in \N$ and bounded by $k$.


\end{abstract}
\section{Introduction}\label{S1}
Reversibility is a property considered in several classes of structures:
a topological space $\X=\la X, \CO\ra$ is reversible iff
$\X \not\cong \la X, \CO '\ra$, for each weaker topology $\CO ' \varsubsetneq \CO$
iff each continuous bijection $\X \rightarrow \X$ is a homeomorphism (or, equivalently, open);
a  relational structure $\X=\la X,\bar \r \ra$ is reversible
iff $\la X,\bar \r \ra \not\cong\la X,\bar \s \ra$,
for each smaller interpretation $\bar \s\varsubsetneq \bar \r$ of the language,
iff every bijective endomorphism (condensation) $f:\X \rightarrow \X$ is an automorphism, namely,
iff $\Cond (\X )=\Aut (\X)$, where $\Cond (\X)$ denotes the set of all self-condensations of $\X$.

It turns out that many relevant structures are reversible
(e.g.\ the spaces $\BR ^n$, linear orders, Boolean algebras \cite{Kuk}, Henson graphs \cite{KuMo2}, Henson digraphs \cite{KRet})
and that, on the other hand, several relevant structures do not have that property
(e.g.\ normed spaces of infinite dimension \cite{Raj}, the random graph, the random poset \cite{KuBfs}).

In this paper we continue the investigation from \cite{KuBfs},
where several classes of non-reversible posets (mainly lattices) were detected
(including homogeneous-universal posets,
the divisibility lattice $\la \N ,\,\mid\,\ra$,
the ideals $[\k ]^{<\l}$,
the ideal of meager Borel subsets of the Baire space,
the direct powers of rationals and integers, $\Q ^\k$ and $\Z ^\k$).
Here we consider the phenomenon of reversibility in the class of trees.

Well orders are simple examples of reversible trees
and by the following result of Kukiela finite disjoint unions of well orders are reversible too.
\begin{fac}[\cite{Kuk}]\label{T102}
Well founded posets with finite levels are reversible.
\end{fac}
The following characterization shows that even for infinite disjoint unions of finite well orders $\bcd _{i\in I}n_i$ the situation is less trivial.
By \cite{KuMo5} a sequence $\la n _i :i\in I\ra \in \N ^I$ is called a {\it reversible sequence of natural numbers}
iff there is no non-injective surjection $f:I\rightarrow I$
such that $\sum _{i\in f^{-1}[\{ j \}]}n_i =n _j$, for all $j\in I$,
iff the tree $\bcd _{i\in I}n_i$ is reversible.
For a characterization of such sequences
we recall that a set $K\subset\N$ is called {\it independent}
iff $K=\emptyset$ or $n\not\in \la K\setminus \{ n \}\ra$, for all $n\in K$,
where $\la K\setminus \{ n \}\ra$ is the subsemigroup of the semigroup  $\la \N , +\ra$
generated by the set $K\setminus \{ n \}$. Let $\gcd (K)$ denote the greatest common divisor of the numbers from $K$.
\begin{fac}[\cite{KuMo5}]\label{TC001}
A sequence $\langle n_i :i\in I\rangle \in \N ^I$ is reversible if and only if the set
$K=\{ m\in \mathbb{N} : |\{ i\in I : n_i=m\}|\geq \o\}$
is independent and, if $K\neq \emptyset$, then $\gcd (K)$ divides at most finitely many elements
of the {\em set} $\{ n_i :i\in I \}$.
\end{fac}
We note that independent sets are finite and that the set $\{ n,n+1, \dots ,2n-1\}$ is independent, for each $n\in \N$.
Regarding the general situation, when $\BL _i$, $i\in I$, are pairwise disjoint well orders,
$\BL _i \cong \a _i =\gamma_i+n_i\in \Ord$, where $\gamma_i\in \Lim \cup \{0\}$ and $n_i\in\omega$, for $i\in I$,
and defining $I_\alpha   :=  \{i\in I:\alpha_i=\alpha\}$, for $\a\in \Ord$, and $J_\gamma   :=  \{j\in I:\gamma_j=\gamma\}$, for $\g\in \Lim \cup \{0\}$, we have
\begin{fac}[\cite{KuMo2}]\label{T000}
$\bcd _{i\in I} \BL_i$ is a reversible tree iff
$|I_\alpha |<\omega$, for all $\alpha \in \Ord$,
or
there exists $\g =\max \{ \g _i : i\in I\}$,
$|I_\a|<\o  $, for all $\a \leq \g$,
and $\la n_i : i\in J_\g \setminus I_\g\ra$ is a reversible sequence of natural numbers.
\end{fac}
For example $(\bcd _{n\in \N}n ) \du \o \du (\bcd _\o (\o +4)) \du (\bcd _\o (\o +6))\du (\bcd _{k\in \o}(\o +2k +1))$ is a reversible tree,
while $\bcd _\mu \g$ is not, if $\g$ is a limit ordinal and $\mu$ an infinite cardinal.

Although the connected (i.e.\ rooted) trees are more frequently considered and constructed in several contexts,
the disconnected ones are important in our analysis.
Namely, by Fact \ref{T8111}, a tree $\BT$ is reversible iff $N\up$ is a reversible tree, for each node $N$ of $\BT$;
and, clearly, the tree $N\up$ is disconnected, whenever $|N|>1$.
Generally speaking, denoting the language with one binary relation by $L_b$, we have
\begin{fac}[\cite{KuMo4}]\label{TB051}
(a) An $L _b$-structure is reversible iff each union of its components is reversible.
In particular, the connectivity components of a reversible $L_b$-structure are reversible (\cite{KuMo4} Theorem 3.1).

(b) An $L_b$-structure with finitely many components is reversible iff all its components are reversible (\cite{KuMo4} Corollary 3.4).
\end{fac}
Thus an important question in our context is: when is an infinite disjoint union $\bcd _{i\in I}\BT _i$ of connected reversible trees reversible?
The following statement (see Corollary 5.2  and Proposition 5.6 of \cite{KuMo4})
provides sufficient conditions for reversibility expressed in terms of cardinal and ordinal invariants of trees.
\begin{fac}[\cite{KuMo4}]\label{T001}
If $\BT _i$, $i\in I$, are pairwise disjoint connected reversible trees
and the double sequence $\la \la|T_i|,\he (\BT_i)\ra:i\in I\ra$ is finite-to-one,
then $\bcd _{i\in I}\BT _i$ is a reversible tree. This holds if, in particular,
the sequence of cardinals $\la |T_i|:i\in I\ra$ is finite-to-one or the sequence of ordinals $\la \he (\BT _i):i\in I\ra$ is finite-to-one.
\end{fac}
So the remaining situation is when there are an infinite set $J\subset I$, an ordinal $\a$ and a cardinal $\k$
such that all the trees $\BT _i$, for $i\in J$, are of height $\a$ and of cardinality $\k$.
Some examples of such unions are ``bad trees" (having all branches of height $\he (\BT)=|T|=|L_0|$, where $|T|$ is a regular cardinal)
and in Section \ref{S3}, using a back and forth characterization of (non-)reversible relational structures,
we show that bad trees are not reversible.
Consequently, if $\BT$ is a countable tree of height $\o$ without leafs, then $\BT$ is reversible
iff all nodes of $\BT$ are finite
iff the set $[T]$ of branches of $\BT$ is a compact subset of the Baire space.

In Section \ref{S4} we show that a tree $\BT$ is non-reversible
iff it contains at least one of two substructures of a specific form.
The first one is a {\it critical node}, $N\in \CN (\BT)$, where  $N\up$ is a tree with infinitely many components and an additional property.
Thus, the results concerning disconnected structures are relevant here.

The second one is called an {\it archetypical subtree of $\BT$};
this is a tree having the last two levels isomorphic to the simplest, ``archetypical", non-reversible tree: $\bcd _\o 1 \du \bcd _\o 2$,
and having some additional properties.

If $\BT$ is a tree with finite nodes, then it does not contain critical nodes, thus it is reversible iff it does not contain archetypical subtrees.
Using that characterization in Section \ref{S5} we prove that
if $\he (\BT)\leq \o +1$ or for each ordinal $\a \in [\o ,\he (\BT ))$
the sequence $\la|N|: \CN (\BT ) \ni N\subset L_\a\ra $ is almost constant (mod finite)
or the sequence $\la\la |N|,|N\up|\ra : \CN (\BT ) \ni N\subset L_\a\ra $ is finite-to-one,
then $\BT$ is reversible.
Applying that result we construct several reversible trees and prove that all regular $n$-ary trees are reversible,
that there are reversible Aronszajn trees,
and that the existence Suslin or Kurepa trees implies the existence of reversible ones.

In Section \ref{S6} we regard disconnected trees
and using the methods described above prove that
the tree $\bcd _\mu {}^{<\a }\l$ is reversible iff $\min\{\a ,\l\mu\} <\o$,
where $\l >1$ and $\mu >0$ are cardinals and $\a>0$ is an ordinal.
Also, for $\l \geq \o$, the tree $\bcd _{i\in I}{}^{<n_i}\l$ is reversible
iff the sequence $\la n_i :i\in I\ra $ is finite-to-one or bounded by some $k\in \N$ and almost constant at $k$ (i.e.\ $|\{ i\in I:n_i\neq k\}|<\o$).
\section{Preliminaries}\label{S2}
\paragraph{Basic definitions and notation}
A partial order $\BT =\la T,<\ra$ is a {\it tree} iff  $( \cdot  ,t):= \{ s\in T :s<t\}$ is a well order or the empty set, for each $t\in T$.
The ordinal $\he (t):=\otp ( ( \cdot  ,t))$ is the {\it height} of the element $t$ of $T$ (here $\otp (\emptyset)$ is the ordinal 0)
and, for an ordinal $\a$, the set $\Lev _\a (\BT)=\{ t\in T : \he (t) =\a\}$ is the {\it $\a$-th level} of $\BT$.
The least ordinal $\a$ satisfying $\Lev _\a (\BT)=\emptyset$ is the {\it height of the tree} $\BT$, in notation $\he (\BT)$.
Instead of $\Lev _\a (\BT)$ we will write $L_\a$, whenever the context admits it.
For $\a \leq \he (\BT)$, by $\BT |\a$ we denote the tree $\bigcup _{\b<\a}L_\b$.
Maximal chains of $\BT$ are the {\it branches} of $\BT$ and the {\it height of a branch} $b$ is the ordinal $\otp (b)$.
Maximal elements of $\BT$ (if they exist) are the {\it leafs of $\BT$} and the least element of $\BT$, if it exists, is the {\it root of $\BT$};
then $\BT$ is a {\it rooted tree}.

The transitive closure $\r$ of the relation $\leq \cup \geq$
is the minimal equivalence relation on the set $T$ containing the relation $<$
and it is called the {\it connectivity relation on $\BT$}.
It is easy to see that $t \mathrel{\r } t'$ iff there is $s\in T$ such that $s\leq t,t'$.
The corresponding equivalence classes are the {\it connectivity components of $\BT$}.
If $|T/\r |=1$, the tree $\BT$ is {\it connected}; this holds iff $\BT$ is a rooted tree, iff $|L_0|=1$.
If $\BT _i$, $i\in I$, are pairwise disjoint rooted trees, $\bcd _{i\in I}\BT _i$ denotes their {\it disjoint union};
$\bcd _I \BT$ denotes the disjoint union of $|I|$-many disjoint copies of $\BT$;
if $\BL$ is a well order, $\BL +\BT$ denotes the corresponding ordinal sum ($L < T$).

If $\sim $ is the equivalence relation on the set $T$ defined by $s\sim t$ iff $(\cdot ,s)=(\cdot ,t)$,
the corresponding equivalence classes are the {\it nodes} of $\BT$.
Clearly, nodes are subsets of levels of $\BT$,
the set of immediate successors $\is (t)$ of a non-maximal element $t$ of $T$ is a node of $\BT$,
but the nodes on limit levels are not of that form.
The set $T/\!\sim$ of nodes of $\BT$ will be denoted by $\CN (\BT)$.

For $t\in T$, $( \cdot  ,t]:=\{ s\in T :s\leq t\}$,
$[t , \cdot ):=\{ s\in T :s\geq t\}$ and  $(t , \cdot ):=\{ s\in T :s> t\}$.
For $S\subset T$ we define $S\up :=\{ t\in T : \exists s\in S \;\; t\geq s\}$ and $S\down :=\{ t\in T : \exists s\in S \;\; t\leq s\}$;
the set $S$ is {\it downwards closed} iff $S \down =S$.
By $\parallel$ we denote the incomparability relation: $s \parallel t$ iff $s\not\leq t$ and  $t\not\leq s$.

For a cardinal $\l >0$, a tree $\BT$ is called {\it $\l$-ary} iff $|N|=\l$, for each node $N$ of successor height
(i.e., iff $|\is (t)|=\l$, for each non-maximal $t\in T$).
$\BT$ is called {\it regular} iff it has a root and $|N|=1$, for each node $N$ of limit height.
It easy to prove that a regular $\leq \l$-ary tree of height $\a $
is isomorphic to a downwards closed subtree of the {\it complete $\l$-ary tree of height $\a$}, ${}^{<\a }\l=\bigcup _{\b <\a }{}^{\b }\l$.
\paragraph{Elementary facts} The following facts will be used throughout the paper.
\begin{fac}\label{T8124}
If $\;\BT $ is a tree and $f\in \Cond (\BT)$, then

(a) $\he (t)\leq \he (f(t))$, for all $t\in T$;

(b) $f\in \Aut (\BT)$ iff  $\;\he (t)= \he (f(t))$, for all $t\in T$;

(c) $T |\a \subset f[T |\a]$, for each $\a \leq \he (\BT)$.
\end{fac}
\dok
(a) Let $\he (t)=\a$ and $\he (f(t))=\b$.
For $s\leq s' <t$ we have $f(s)\leq f(s')<f(t)$;
thus $f[( \cdot  ,t)]\subset ( \cdot  ,f(t))$
and $f\upharpoonright ( \cdot  ,t)$ is a monotone injection from the well order $( \cdot  ,t)$ into the well order $(\cdot  ,f(t))$
and, hence, an embedding,
which implies that $\a \leq \b$.

(b) If $f\in \Aut (\BT)$, then, as above, $f\upharpoonright ( \cdot  ,t):( \cdot  ,t)\rightarrow (\cdot  ,f(t))$ is an embedding.
Moreover, if $u<f(t)$, then $u=f(s)$, for some $u\in T$, and, since $f\in \Aut (\BT)$, $s<t$.
Thus  $f[( \cdot  ,t)]=(\cdot  ,f(t))$ and, hence, $( \cdot  ,t)\cong (\cdot  ,f(t))$ which gives $\he (t)= \he (f(t))$.

If $f\not\in \Aut (\BT)$ there are $s,t\in T$ such that $f(s)<f(t)$ and $s\not<t$.
Since $t<s$ would imply that $f(t)<f(s)$ we have $s\| t$.
If $\he (s)< \he (f(s))$ or $\he (t)< \he (f(t))$ we are done.
Otherwise by (a) we have $\he (s)= \he (f(s))< \he (t)= \he (f(t))$
and, since $s\| t$, there is $u<t$ such that $\he (u)=\he (s)$ and $u\| s$,
which implies $f(u)\neq f(s)$.
Since $u<t$ we have $f(u)<f(t)$
and, since  $f(s)<f(t)$ as well,
$f(u)$ and $f(s)$ are different comparable elements of $T$.
So, since $\he (f(u))\geq \he (u)=\he (s)=\he (f(s))$ we have $\he (f(u))> \he (u)$.

(c) Let $t\in T |\a$.
Since $f$ is a surjection there is $s\in T$ such that $t=f(s)$
and by (a) we have $\he (s)\leq \he (f(s))=\he (t) <\a$.
Thus $s\in T |\a$ and, hence, $t=f(s)\in f[T |\a]$.
\hfill $\Box$
\begin{fac}\label{T8130}
If $\BT$ is a non-reversible tree
and if for $f \in \Cond (\BT )\setminus \Aut (\BT )$ we define $\a _f :=\min \{ \he (t): t\in T \land \he (t)< \he (f(t))\}$, then
\begin{eqnarray}\label{EQ8157}
\a _\BT & := & \min \Big\{ \a _f: f \in \Cond (\BT )\setminus \Aut (\BT )\Big\}\\
        & =  & \max \Big\{ \b<\he (\BT) : \forall f \in \Cond (\BT ) \;\; f\upharpoonright (\BT |\b)\in \Aut (\BT |\b )\Big\}\nonumber.
\end{eqnarray}
\end{fac}
\dok
First we show that $\a _\BT \in S:= \{ \b<\he (\BT) : \forall f \in \Cond (\BT ) \; f\upharpoonright (\BT |\b)\in \Aut (\BT |\b )\}$.
Clearly, $\a _\BT <\he (\BT)$.
If $f \in \Cond (\BT )$ and $t\in T |\a _\BT$,
then $\he (t) <\a _\BT \leq \a _f$
and, hence, $\he (t) = \he (f(t))$
which implies that $f(t)\in T |\a _\BT$.
Thus $f[ T |\a _\BT]\subset T |\a _\BT$
and, by Fact \ref{T8124}(c), $f[ T |\a _\BT]= T |\a _\BT$.
Since $\he (t) = \he (f(t))$, for all  $t\in T |\a _\BT$,
by Fact \ref{T8124}(b) we have $f\upharpoonright (\BT |\a _\BT)\in \Aut (\BT |\a _\BT )$;
thus, $\a _\BT \in S$.

Second, for $\b >\a _\BT$ we prove that $\b \not\in S$.
Let $f \in \Cond (\BT )\setminus \Aut (\BT )$, where  $\a _f =\a _\BT$,
and let $t\in L_{\a _\BT}$ be such that $\a _\BT =\he (t) < \he (f(t))$.
Since $\b >\a _\BT$ we have $t\in T|\b$.
If $f[T|\b ]\not\subset T|\b$, then, clearly, $f\upharpoonright (\BT |\b)\not\in \Aut (\BT |\b )$.
Otherwise, by Fact \ref{T8124}(c) we have $f[ T |\b ]= T |\b$
and since $t\in T|\b$ and $\he (t) < \he (f(t))$ by Fact \ref{T8124}(b) we have $f\upharpoonright (\BT |\b)\not\in \Aut (\BT |\b )$ again.
Thus, $\b \not\in S$ indeed.
\hfill $\Box$
\begin{rem}\label{R8105}\rm
If $\BT$ is a non-reversible tree, then using Fact \ref{T8130} it is easy to show that $|L_{\a _\BT}|\geq \o$.
We omit a proof here since this follows from Proposition \ref{T8108}.

For each ordinal $\a$ there is a simple tree $\BT$ such that $\a _\BT=\a$: take $\BT =\a +(\bcd _\o 1 \du \bcd _\o 2) $.
Less trivial examples are described in Example \ref{EX8104}.

By Fact \ref{T8130}, if $\BT$ is a non-reversible tree, then $f\upharpoonright (\BT |\a _\BT)\in \Aut (\BT |\a _\BT )$, for each $f \in \Cond (\BT )$.
But Example \ref{EX8104} also shows that it is possible that the restriction $\BT |\a _\BT$ is not reversible.
This means that there is a condensation $g\in \Cond (\BT |\a _\BT)\setminus \Aut (\BT |\a _\BT )$ which has no extension to a condensation of $\BT$.
\end{rem}
\begin{fac}\label{T8111}
If $\BT$ is a tree, then we have

(a) $\BT$ is reversible $\Leftrightarrow$ $\forall N\in \CN (\BT ) \;\forall S\in P(N)\setminus \{ \emptyset \} \;( S\up \mbox{is a reversible tree})$;
    then, in particular,  $t\up$ is a reversible tree, for all $t\in T$;

(b) If $|L_0|<\o$, then we have: $\BT$ is reversible $\Leftrightarrow$ $\forall t\in T \,(t\up \mbox{is a reversible tree})$.
\end{fac}
\dok
(a) The implication ``$\Leftarrow$" is trivial, because $L_0\in \CN (\BT )$ and $L_0 \up =T$.

Conversely, assuming that $S\subset N\in \CN (\BT )$ and $f\in \Cond (S\up)\setminus \Aut (S\up)$, we prove that
$g:= \id _{T\setminus S\uparrow }\cup f\in \Cond (\BT )\setminus \Aut (\BT)$.
It is clear that $g:T \rightarrow T$ is a bijection. For $x,y\in T$ satisfying $x<y$ we show that $g(x)<g(y)$.
This is evident if $x,y\in S\up$ or $x,y \in T\setminus S\up $
and, since that set $S \up$ is upwards closed, $x\in S\up \not\ni y$ is impossible.
So, if $x\not\in S\up \ni y$,
then there is $t\in S$ such that $t\leq y$
and $x<y$ implies that $x<t$.
Now  $f\in \Cond (S\up)$ implies that $f(y)\in S\up$.
Let $t_1 \in S$, where $t_1 \leq f(y)$.
Since $x\in (\cdot ,t)=(\cdot ,t_1)$ we have $g(x)=x <t_1 \leq f(y)=g(y)$.
Thus we have $g\in \Cond (\BT )$.

Since $f\not\in \Aut (S\up )$ there are $x,y\in S\up$ such that $f(x)<f(y)$
and, hence, $g(x)<g(y)$, but $x\not< y$. So $g\in \Cond (\BT )\setminus \Aut (\BT )$ and $\BT$ is not reversible.

(b) The implication ``$\Rightarrow$" follows from (a). If $t\up$ is a reversible tree, for each $t\in T$,
then, since $T =\bigcup _{t\in L_0} t\up$ is a finite union of pairwise disjoint, connected and reversible structures,
by Fact \ref{TB051}(b) the tree $\BT$ is reversible too.
\hfill $\Box$
\begin{ex}\label{EX8100}\rm
The assumption $|L_0|<\o$ in Fact \ref{T8111}(b) can not be omitted.
If $\BT$ is $\bcd _\o 1 \du \bcd _\o 2$ or any non-reversible (thus infinite) disjoint union of well orders (see Fact \ref{T000}),
then $t\up$ is a reversible tree (in fact, a well order), for each $t\in T$.

If a tree $\BT$ is not reversible, $t\in T$ and the principal filter $t\up =[t,\cdot )$ is not reversible,
then $[s,\cdot )$ is not reversible, for each $s\leq t$. Thus $\I _\BT := \{ t\in T : [t,\cdot )$ is not reversible$\}$
is a downwards-closed subset of $T$.
If $\BT$ is as above, we have $\I _\BT =\emptyset$.
On the other hand, if $T={}^{<\o}\o$, then $\I _\BT =T$, see Example \ref{EX8103}.
\end{ex}
The following fact shows how the statements and examples concerning reversibility of disconnected trees
can be converted into the results about connected (i.e.\ rooted) trees.
If $\BT=\la T,<\ra$ is a tree, $r\not\in T$ and $T_r:= \{ r\} \cup T$, by $\BT _r$ we denote the rooted tree $\la T_r ,(\{ r\} \times T) \;\cup <\ra$ and we have
\begin{fac}\label{T8100}
A disconnected tree $\BT$ is reversible iff the tree $\BT _r$ is reversible.
\end{fac}
\dok
Let the tree $\BT _r$ be reversible. Since $L_0 (\BT)=\is (r)\in \CN (\BT _r)$, by Fact \ref{T8111}(a)
the tree $\BT =L_0 (\BT)\up$ is reversible.

Let the tree $\BT _r$ be non-reversible,
$F\in \Cond (\BT _r)$ and $t_0\in T_r$, where $\he _{\BT _r}(t_0)<\he _{\BT _r}(F(t_0))$.
Then, clearly, $F(r)=r$, $f:= F|T \in \Cond (\BT)$ and $t_0\in T$.
So, for a proof that $\BT$ is not reversible it remains to be shown that $\he _{\BT }(t_0)<\he _{\BT }(f(t_0))$.
Generally, for $t\in T$ we have
\begin{equation}\label{EQ8168}
\he _{\BT _r}(t) = \left\{
                     \begin{array}{ll}
                      \he _{\BT }(t) +1 , & \mbox { if }  \he _{\BT }(t)<\o, \\
                      \he _{\BT }(t),      & \mbox { if }  \he _{\BT }(t)\geq\o .
                     \end{array}
                   \right.
\end{equation}
So, if $\he _{\BT _r}(t_0)\geq \o$, then $\he _{\BT }(t_0)=\he _{\BT _r}(t_0)<\he _{\BT _r}(F(t_0))=\he _{\BT }(f(t_0))$.
If $\he _{\BT _r}(F(t_0))< \o$, then $\he _{\BT }(t_0)+1=\he _{\BT _r}(t_0)<\he _{\BT _r}(F(t_0))=\he _{\BT }(f(t_0))+1<\o$,
which gives $\he _{\BT }(t_0)<\he _{\BT }(f(t_0))$.
Finally, if $\he _{\BT _r}(t_0)< \o$ and $\he _{\BT _r}(F(t_0))\geq \o$,
then we have $\he _{\BT }(t_0)+1=\he _{\BT _r}(t_0)<\o \leq \he _{\BT _r}(F(t_0))=\he _{\BT }(f(t_0))$.
\hfill $\Box$
\section{Trees with all branches of maximal height. Bad trees}\label{S3}
In this section we detect a class of non-reversible trees
using a characterization of non-reversible structures from \cite{KuBfs}.
First we recall that a partial order $\P =\la P, \leq \ra $ is called {\it $\k$-closed} (where $\k\geq \o$ is a cardinal)
iff whenever $\g <\k$ is an ordinal and $\la p_\a : \a <\g \ra$ is a sequence in $P$ such that $p_\b \leq p_\a$, for $\a <\b<\g$,
there is $p\in P$ such that $p\leq p_\a$, for all $\a <\g$.

Let $L=\la R_i:i\in I\ra$ be a relational language, where $\ar (R_i)=n_i\in \N$, for $i\in I$,
and let $\X $ be an $L$-structure.
A function $\f$ is called a {\it partial condensation of $\X$}, we will write $\f \in \Pc (\X )$,
iff $\f$ is a bijection which maps $\dom (\f )\subset X$ onto $\ran (\f )\subset X$ and
\begin{equation}\label{EQ0092}
\forall i\in I \;\; \forall \bar x \in (\dom (\f ))^{n_i} \;\;(\bar x \in R_i^\X \Rightarrow \f\bar x \in R_i^\X ).
\end{equation}
A non-empty set $\Pi \subset \Pc (\X )$ is called a {\it back and forth system of condensations} (shortly: b.f.s.) iff

(bf1) $\forall \f\in \Pi \;\; \forall a\in X \;\;  \exists \p \in \Pi\;\;(\f \subset \p \land a\in \dom \p )$,

(bf2) $\forall \f\in \Pi \;\; \forall b\in X \;\;\,\exists \p \in \Pi\;\;(\f \subset \p \land b\in \ran \p )$.

\noindent
We will say that a b.f.s.\ $\Pi $ is a {\it $\k$-closed b.f.s.}
iff the partial order $\la \Pi ,\supset\ra$ is $\k$-closed.
A partial condensation $\f ^* \in \Pc (\X )$ is called {\it bad}
iff  $\f ^*$ is not a partial isomorphism of $\X$;
that is, there are $i\in I$ and $\bar x =\la x_0 ,\dots ,x_{n_i-1}\ra\in \dom(\f )^{n_i}$
such that $\bar x \not\in R_i^\X$ and $\f ^*\bar x \in R_i^\X$.
\begin{fac}[\cite{KuBfs}]\label{T8196}
A relational structure $\X $ of size $\k\geq \o$ is not reversible if and only if there is a $\k$-closed b.f.s.\ $\Pi \subset \Pc (\X )$ containing a bad condensation.
\end{fac}
A tree $\BT$ will be called {\it bad} iff

(b1) $\he (\BT ) =|T|=|L_0|$ is a regular cardinal, and

(b2) each branch of $\BT$ is of height $\he (\BT )$.
\begin{te}\label{T8110}
(a) Bad trees are not reversible.

(b) If  $\BT$ is a tree having a node $N$ such that $S\up$ is a bad tree, for some $S\subset N$, then $\BT$ is not reversible.
\end{te}
\dok
(a) Let $\BT$ be a bad tree of size $\k$.
We show that the assumptions of Fact \ref{T8196} are satisfied.
Let $a_0,a_1\in L_0$, where $a_0\neq a_1$. Since $a_0\not\in \Max (T)$ there is  $b_0\in \Min ( (a_0 , \cdot))$.
Then $\f ^*:=\{ \la a_0, b_0\ra ,\la a_1, a_0\ra\}\in \Pc (\BT )$ is a bad condensation and $\dom (\f ^*)=\{ a_0, a_1\}$
is a downwards closed set; thus $\f ^* \in \Pi$, where
\begin{equation}\label{EQ8117''}
\Pi :=\{ \f \in \Pc (\BT ) : |\f |<\k \land (\dom \f ) \!\downarrow \;= \dom \f \land \f ^* \subset \f \} .
\end{equation}
Since $\k$ is a regular cardinal and the union of  downwards closed sets is downwards closed, the poset $\la \Pi ,\supset\ra$ is $\k$-closed.

(bf1) Let $\f \in \Pi $ and $a\in T \setminus  \dom \f$.

{\it Case 1}: $(\cdot , a]\cap \dom \f\neq \emptyset$.
Then the set $D=\dom \f \cup (\cdot ,a]$ is downwards closed,
$a\in (\cdot , a]\setminus \dom \f$ and defining
\begin{equation}\label{EQ8126''}
u := \min \Big( (\cdot , a]\setminus \dom \f \Big),
\end{equation}
we have $D= \dom \f \;\dot{\cup}\; [u,a]$.
Since the set $\dom \f$ is downwards closed
we have $(\cdot ,u)\subset \dom \f$
and, hence, $\f [( \cdot ,u)]$ is a chain in $\ran \f$.
If $C$ is a maximal chain in $\ran \f$ such that $\f [(\cdot ,u)]\subset C$,
then, since $|\ran \f |=|\f |<\k$, we have $|C|<\k$.
Let $B$ be a branch in $\BT$ containing $C$.
Then by (b2) we have $\otp (B)=\k$
and, hence, there is $z\in B$ such that $C<z$.
Since $\otp ([u, a])<\k$ there is an embedding $\eta :[u, a]\hookrightarrow B$ such that $\eta (u)=z$.
Let $\p :D\rightarrow T$, where $\p \upharpoonright \dom \f =\f$ and $\p \upharpoonright [u,a]=\eta$.
It is evident that $\p$ is an injection.

In order to prove that $\p$ is a homomorphism,
assuming that $x,y\in D$ and $x< y$ we show that $\p (x)< \p (y)$.
Since the set $\dom \f$ is downwards closed it is impossible that $x\not\in \dom \f \ni y$.
If $x,y\in \dom \f $, then $\p (x)=\f (x)< \f (y)=\p (y)$ ($\f$ is a homomorphism).
If $x,y\not\in \dom \f $, then $x, y\in [u,a]$
and, since $\eta$ is an embedding, $\p (x)=\eta (x)< \eta (y)=\p (y)$.
If $x\in \dom \f \not\ni y$,
then $y\in [u,a]$ and, hence, $x\in (\cdot , a]\cap \dom \f=(\cdot , u)$.
So, $\f (x)\in \f[(\cdot , u)]$
and $\p (x)=\f (x)<z=\eta (u)\leq \eta (y)=\p (y)$. Thus $\f \subset \p \in \Pi$ and $a\in \dom \p$.

{\it Case 2}: $(\cdot , a]\cap \dom \f= \emptyset$. Then again $D=\dom \f \cup (\cdot ,a]$ is a downwards closed set,
and, hence, defining $u =\min (\cdot , a]$, we have $[u,a]\parallel \dom \f$.
Let $\he (a)=\a$.
Since $|\ran \f|< \k$ there is $z\in L_0$ such that
$[z ,\cdot) \cap \ran \f =\emptyset$ and, by (b2) there is $\eta :[u, a]\hookrightarrow [z ,\cdot)$ such that $\eta (u)=z$.
Now $\p =\f \cup \eta :D\rightarrow T$ is an injection, $a\in \dom \p$ and as above we show that it is a homomorphism.

(bf2) Suppose that $\f \in \Pi$ and $b\in T \setminus  \ran \f$.
Since $|\dom \f |<\k$ and $|L_0|=\k$ there is $a\in L_0 \setminus \dom \f$.
Then $\p := \f \cup \{ \la a ,b \ra\}$ is an injection and,
since $a\in L_0$, the set $\dom \p =\{ a\} \cup \dom \f$ is downwards closed.

For a proof that $\p$ is a homomorphism,
assuming that $x,y\in \{ a\} \cup \dom \f$ and $x< y$ we show that $\p (x)< \p (y)$.
Since the set $\dom \f$ is downwards closed and $\dom \p \setminus \dom \f =\{ a\}$,
it is impossible that $x\not\in \dom \f \ni y$ or $x,y\not\in \dom \f $.
$x\in \dom \f \not\ni y$ would imply that
$y=a$ and, hence $x<a$, which is impossible because $a\in L_0$.
Thus $x,y\in \dom \f $ and, again, $\p (x)=\f (x)< \f (y)=\p (y)$.
So, $\f \subset \p \in \Pi$ and $a\in \dom \p$.

(b) If $S\subset N\in \CN (\BT)$ and $S\up$ is a bad tree,
then by (a) $S\up$ is not reversible
and, by Fact \ref{T8111}(a), $\BT$ is not reversible.
\hfill $\Box$
\begin{rem}\label{R8103}\rm
By definition, a bad tree of size $\k$ is a disjoint union of the form $\bcd _{i<\k}\BT _i$,
where $\BT _i$ is a connected tree and $|T_i|=\he (\BT _i)=\k$, for each $i\in \k$.
By Fact \ref{TB051}, Theorem \ref{T8110} says something new if the trees $\BT _i$, $i<\k$, are reversible.
The simplest case is when $\BT _i \cong \k$, for all $i\in\k$.
Also, by Theorem \ref{T8119}(a) the tree ${}^{<\k }2$ is reversible and if $2^{<\k}=\k$ (in particular, under the GCH) it is of size $\k$;
thus, the tree $\bcd _\k {}^{<\k }2$ is not reversible by Theorem \ref{T8110}(a).
But if, for example, $2^\o=\o _2$, then $|{}^{<\o _1 }2|=\o _2$ and the tree $\bcd _{\o_1} {}^{<\o_1 }2$ is not bad.
We note that, by Theorem \ref{T8126}, all trees of the form $\bcd _\mu {}^{<\a }\l$, where $\a, \mu \geq \o$,  are non-reversible.
\end{rem}
\begin{te}\label{T8197}
Let $\BT$ be a tree with all branches of height $\he (\BT)$.

(a) If $\he (\BT )<\o$, then $\BT$ is reversible;

(b) If $\he (\BT)=|T|=\o$, then $\BT$ is reversible iff all nodes of $\BT$ are finite
iff, regarding $\BT$ as a downwards closed subtree (initial part)
of the tree ${}^{<\o}\o$, the set $[T]:=\{ x\in {}^{\o }\o :\forall n\in \o \; x\upharpoonright n \in T\}$ is a compact subset of the Baire space $\o ^\o$.
\end{te}
\dok
(a) Let $\he (\BT) =n \in \N$ and $\he (b)=n$, for each branch $b$ of $\BT$.
Suppose that the tree $\BT$ is not reversible and $f\in \Cond (\BT )\setminus \Aut (\BT )$. Then by Fact \ref{T8124}(b)
there is $t\in T$ such that $k:=\he (t) < \he (f(t))< n$. If $b$ is a branch of $\BT$ such that $t\in b$, then since $\he (b)=n$ we have
$b\cap [t,\cdot )=\{ t ,t_{k+1}, \dots ,t_{n-1}\}$, where $t_i\in b\cap L_i$,
and  $t <t_{k+1}< \dots <t_{n-1}$, which implies that $f(t) <f(t_{k+1})< \dots <f(t_{n-1})$
and, hence, $k<\he (f(t)) <\he (f(t_{k+1}))< \dots <\he (f(t_{n-1}))$.
Thus $\he (f(t_{n-1}))\geq n$, which is false.

(b) If $N$ is a node of $\BT$ and $|N|=\o$, then  $N\up$ is a bad tree of size $\o$ and, by Theorem \ref{T8110}(b), $\BT$ is not reversible.

Let $|N|<\o$, for each $N\in \CN (\BT)$. Since $L_0$ is a node of $\BT$ we have $|L_0|<\o$. If $|L_n|<\o$, then, since $L_{n+1}=\bigcup _{x\in L_n}N_x$, where
$N_x$ is the set of immediate successors of $x$, we have $|L_{n+1}|<\o$. Thus $\BT$ is reversible by Fact \ref{T102}.

By Fact \ref{T8100} we can suppose that $\BT$ is a rooted tree. A simple induction shows that, under the assumptions,
$\BT$ is isomorphic to (and, hence, can be identified with) a downwards closed subtree of ${}^{<\o}\o$ with no leafs (a ``pruned tree").
It is known (see \cite{Kech}) that the mapping $\BT \mapsto [T]$ is a bijection
from the set $\Prt$ of pruned subtrees of ${}^{<\o}\o$
onto the set $\F$ of closed subsets of the Baire space
and its inverse is given by $F\mapsto \Tree (F):=\{ x\upharpoonright n :x\in F \land n\in \o\}$.
In addition, the set $[T]$ is compact iff all nodes of $\BT$ are finite (see \cite{Bart}, p.\ 6),
which proves the second equivalence.
\hfill $\Box$
\begin{ex}\label{EX8102}\rm
By  Theorem \ref{T8197}(a), the tree $\bcd _{i\in I}{}^{<n}\l_i$ is reversible, whenever $n\in\N$, $I\neq \emptyset$ and $\la \l _i:i\in I\ra$ is a $I$-sequence of non-zero cardinals.
So, the trees ${}^{<n}\o$, $n\in \N$, are reversible, but ${}^{<\o}\o$ is not reversible by Theorem \ref{T8197}(b).

The condition that all branches are of height $\he (\BT)$ can not be removed from Theorem \ref{T8197}(a):
the tree $\bcd _\o 1 \du \bcd _\o 2$ is not reversible.
\end{ex}
\section{Critical nodes and archetypical subtrees}\label{S4}
In this section we give a combinatorial characterization of reversibility in the class of all trees.
We start with some definitions.
A node $N$ of  a tree $\BT$ will be called {\it critical} iff
\begin{equation}\label{EQ8158}
\exists g\in \Cond (N\up ) \;\; \exists a\in N \;\;\he (a)<\he (g(a)).
\end{equation}
A tree $\BT$ will be called {\it archetypical} iff $T =\bigcup _{n\in \Z} (P _n \cup T _n)$,
where

(a1) $\he (\BT )=\a +2$, where $\a \in \Ord$ and $\a >0$,

(a2) $P_n$, $n\in \Z$, are different branches in $\BT |\a$ of height $\a$,

(a3) $P_n < T_n=\{ a_n,b_n\}$, for all $n\in \Z$,

(a4) $a_n \parallel b_n$, for $n\leq 0$, and $a_n <b_n$, for $n> 0$,

(a5) There is $g\in \Aut (\BT |\a )$, such that $g[P_n]=P_{n+1}$, for all $n\in \Z$.

\noindent
If $\BT$ is a suborder of a tree $\BT _1$
we will say that $\BT$ is an {\it archetypical subtree of $\BT _1$}
iff $(\cdot ,a_0)_{\BT _1}=(\cdot ,b_0)_{\BT _1}=P_0$
and the mapping $F_g:T\rightarrow T$ defined by
\begin{equation}\label{EQ000}
F_g:= g\cup \{ \la a_n,a_{n+1}\ra :n\in \Z \}\cup \{ \la b_n,b_{n+1}\ra:n\in \Z \}
\end{equation}
extends to a condensation of $\BT_1$.
\begin{ex}\label{EX8101}\rm
Let $\BT _n =\la \{ a_n, b_n\}, <_n \ra$, $n\in \Z$, be pairwise disjoint trees such that (a4) holds,
let $\a>0$ be an ordinal and let $\P_n$, $n\in \Z$, be pairwise disjoint copies of $\a$.
The disjoint union $\bcd _{n\in \Z}(\P_n +\BT _n)$ will be denoted by
$\D _\a =\bcd _{n\in \Z}(\a +\BT _n)$.
Clearly, $\D _\a$ is a disconnected archetypical tree (irregular, since $|L_0|=\o$)
and adding a root to $\D _\a$ produces a rooted (thus connected) archetypical tree.
\end{ex}
\begin{te}\label{T8198}
A tree is non-reversible iff it contains a critical node or an archetypical subtree.
\end{te}
The theorem follows from the following two propositions.
\begin{prop}\label{T8117}
If  $\BT$ is a tree, then we have

(a) If $N$ is a critical node of $\BT$, then $|N|\geq \o$ and $\BT$ is not reversible.

(b) If $\BT$ has an archetypical subtree, then  $\BT$ is  not reversible.
\end{prop}
\dok
(a) Let $g\in \Cond (N\up )$ and $a\in N$, where $\he (a)<\he (g(a))$.
Suppose that $|N|=n\in \N$, say, $N=\{ v_0, \dots ,v_{n-1}\}$.
Since $g\in \Sym (N \up)$, for each $k<n$ there is $u_k\in N\up$ such that $f(u_k)=v_k$.
By Fact \ref{T8124}(a), in the tree $N\up$ we have $\he (u_k) \leq \he (v_k)=0$,
which implies that $u_k\in N$ and, hence, $g[N]=N$.
But this is impossible, because $a\in N$ and $\he (a)<\he (g(a))$.
By Fact \ref{T8111}(a) $\BT$ is not reversible.

(b) Let $T ' =\bigcup _{n\in \Z} P _n \cup \{ a_n ,b_n\}$ be an archetypical subtree of $\BT$
where $P_n$, $a_n$, $b_n$, for $n\in \Z$, $g$ and $F_g$ are as in the definition of archetypical subtree.
Then there is  $f\in \Cond (\BT )$ such that $F_g \subset f$ and we have $(\cdot ,a_0)_{\BT }=(\cdot ,b_0)_{\BT }$, which gives $\he _{\BT}(a_0)=\he _{\BT}(b_0)$;
thus, since $a_1=F_g (a_0)=f(a_0)$,  $b_1=f(b_0)$, $a_1 <b_1$ and $f\in \Cond (\BT )$, by Fact \ref{T8124}(a) we have
$$
\he _{\BT}(b_0)=\he _{\BT}(a_0) \leq \he _{\BT}(a_1)<\he _{\BT}(b_1)=\he _{\BT}(f(b_0)).
$$
So, by Fact \ref{T8124}(b), $f\not\in \Aut (\BT )$ and the tree $\BT$ is  not reversible.
\hfill $\Box$
\begin{prop}\label{T8108}
Let $\BT$ be a non-reversible tree, $f\in \Cond (\BT)\setminus \Aut (\BT)$
and  $\a  :=\min \{ \he (t): t\in T \land \he (t)< \he (f(t))\}$.
If $\a =0$, then $L_0$ is a critical node of $\,\BT$. Otherwise, we have:

(a) There are $a_0,b_0\in L_{\a }$ such that $f(b_0)> f(a_0) \in L_{\a }$;

(b) $f\upharpoonright (T|\a )\in \Aut (\BT |\a )$;

(c) $f \upharpoonright (\cdot , a_0) : (\cdot , a_0) \rightarrow (\cdot , f(a_0))$ is an isomorphism;

(d) $a_0$ and $b_0$ are in the same node $N_0$ of $\BT$;

(e) The sets $P_n:= f^n[(\cdot , a_0)]$, $n\in \Z$, are branches in the tree $\BT|\a $ of height $\a $;

(f) If $S_n:=\{ s\in T : P_n <s \}$,  then $f|S_n \in \Cond (S_n ,S_{n+1})$, for each $n\in \Z$;

(g) $f^n |S_0\in \Cond (S_0 ,S_n )$, for each $n\in \o$;

(h) If $P_n= P_0$, for some $n\in \N$, then $N_0$ is a critical node of $\BT$;

(i) If $P_n\neq P_0$, for each $n\in \N$, then $P_m \neq P_n$, for different $m,n\in \Z$, and if
$$
a_n=f^n (a_0), \; b_n= f^n(b_0) \;\mbox{ and }\;  T_n=\{ a_n, b_n\}, \mbox{ for } n\in \Z,
$$
then $\BT _f := \bigcup _{n\in \Z} P_n \cup T_n$ is an archetypical subtree of $\BT$ and  $\he (\BT_f)=\a  +2$;

(j) For each $n\in \Z$ we have $N_n := S_n \cap L_\a =\{ s\in T : (\cdot ,s)=P_n \}\in \CN (\BT )$, $S_n =N_n \up$, $N_{n+1}\subset f[N_n]$ and $|N_{n+1}|\leq |N_n|$;

(k) If $|N_0|<\o$, then $|N_1|<|N_0|$ and there is $n_0\in \N$ such that for each $n\geq n_0$ we have $|N_n|=|N_{n_0}|$ and $|N_n \up|=|N_{n_0} \up|$.
\end{prop}
\dok
(a) Let $b_0 \in L_{\a }$, where $\he (f(b_0))> \a $,
let $c\in L_{\a }\cap (\cdot , f(b_0))$ and $a_0=f^{-1}(c)$.
By Fact \ref{T8124}(a) we have $\he (a_0)\leq \he (f(a_0))=\he (c)=\a $
and $\he (a_0) <\a $ is impossible, by the minimality of $\a $.
Thus $a_0\in L_{\a }$ and $f(a_0)=c < f(b_0)$.

(b) By the minimality of $\a $ and Fact \ref{T8124}(a),
for $\xi <\a $ we have $f[L_\xi]= L_\xi$
and, hence, $f[T|\a ]=T|\a $,
which implies that $f\upharpoonright (T|\a )\in \Cond (\BT |\a )$.
Since $\he (t)=\he (f(t))$, for all $t\in T|\a $,
by Fact \ref{T8124}(b) we have $f\upharpoonright (T|\a )\in \Aut (\BT |\a )$.

(c) If $u\in (\cdot , a_0)$, that is $u<a_0$,
then $f(u)< f(a_0)$;
thus $f[ (\cdot , a_0)]\subset (\cdot , f(a_0))$.
If $v\in (\cdot , f(a_0))$,
then $\xi  :=\he (v)< \he (f(a_0))=\a =\he (a_0)$
and for $u\in (\cdot , a_0)\cap L_\xi $ by (b) and Fact \ref{T8124}(b) we have $\he (f(u))=\xi $
and, hence, $f(u)\in (\cdot , f(a_0))\cap L_\xi  =\{ v\}$, that is $f(u)=v$.
So $f[ (\cdot , a_0)]= (\cdot , f(a_0))$.
By (b) and since the surjective restriction of an isomorphism is an isomorphism too, the claim is true.

(d) Clearly we have $(\cdot , a_0)\cong (\cdot , b_0)\cong \a$.
Suppose that $(\cdot , a_0)\neq (\cdot , b_0)$ and let, say, $u\in (\cdot , a_0)\setminus (\cdot , b_0)$.
Then $\b :=\he (u)<\a$ and for $v\in (\cdot , b_0) \cap L_\b$ we have $u\neq v$ and $u\parallel v$.
Since $u<a_0$ and $v<b_0$ we have  $f(u)<f(a_0)<f(b_0)$ and $f(v)<f(b_0)$,
which implies that $f(u)$ and $f(v)$ are comparable.
So, since $f(u)\neq f(v)$, we have $\he (f (u))\neq \he (f (v))$.
But by (b) and Fact \ref{T8124}(b) we have  $\he (f (u))=\he (u)=\he (v)=\he (f (v))$
and we obtain a contradiction.

(e) By (b) we have $f^n \upharpoonright (\BT |  \a  )= (f\upharpoonright (\BT |  \a  ))^n\in \Aut (\BT | \a  )$.
So, since $(\cdot ,a_0)\subset T| \a  $, we have $f^n[(\cdot , a_0)]\subset T| \a  $
and $f^n[(\cdot , a_0)]=(f^n \upharpoonright (\BT |  \a  ))[(\cdot , a_0)]\cong (\cdot , a_0)\cong  \a  $.
Let $C$ be a maximal chain in $\BT | \a  $ such that $f^n[(\cdot , a_0)]\subset C$.
Then $ \a   \hookrightarrow C \subset \BT | \a  $ implies that $C\cong  \a  $.
If $v\in C$, then $\xi  := \he (v) < \a  $, $C\cap L_\xi  =\{ v \}$
and for $u\in (\cdot ,a_0 )\cap L_\xi $ by Fact \ref{T8124}(b) we have $\xi  = \he (u) = \he (f^n (u))$,
which implies that $ f^n (u)\in f^n[(\cdot , a_0)]\cap L_\xi \subset C\cap L_\xi  =\{ v \}$.
So $v=f^n (u)\in f^n[(\cdot , a_0)]=P_n$
and, thus, $P_n = C$ and $P_n$ is a maximal chain in $\BT| \a  $ of height $ \a  $.

(f) Since the surjective restriction of a condensation is a condensation, we have to prove that
\begin{equation}\label{EQ8160'}
f[S_n ] =S_{n+1}.
\end{equation}
Let $u\in S_n $.
If  $r\in P_{n+1}$ and $t\in P_n$, where $r=f(t)$,
then $t<u$ and, hence, $r=f(t)<f(u)$.
Thus $f(u)>r$, for all $r\in P_{n+1}$,
which means that $f(u)\in S_{n+1}$
and ``$\subset$" in (\ref{EQ8160'}) is true.

Suppose that $v\in S_{n+1} \setminus f[S_n]$ and let $v=f(u)$.
Then $u\not \in S_n$
and, by (e), $\he (v)\geq  \a  $,
which by (b) implies that $\he (u)\geq  \a  $.
Since $u\not \in S_n$ there is $t\in P_n$ such that $t\not< u $
and, since  $u\not< t $, we have $u\parallel t$.
Thus, for $w\in (\cdot ,u] \cap L_{\a } $ we have $w\parallel t$ (otherwise we would have $t<w\leq u$)
and $P_w:=(\cdot ,w)$ is a maximal chain in $\BT | \a  $ of height $ \a  $.
For $\xi  := \otp ((\cdot ,w]\cap (\cdot ,t])$ we have $0\leq \xi  <\he (t) < \a  $
and, hence, $\xi  +1 < \he (t) +1 \leq\a  $.
Taking $x\in P_w \cap L_{\xi  +1}$ and $y\in P_n \cap L_{\xi  +1}$ we have $x\parallel y$.
Since $x<w\leq u$ we have $f(x)< f(w)\leq f(u)=v$
and, since $w\in  L_{\a }$, by Fact \ref{T8124} we have $\he (f(w))\geq  \a  $.
Since $v\in S_{n+1}$ it follows that $(\cdot ,v]\cap (T| \a  )=P_{n+1}\subset (\cdot ,f(w)]$.
By (b) we have $f[P_w] \subset (\cdot , f(w)) \cap (T| \a   )=P_{n+1}$;
since $x\in P_w$ we have $f(x)\in P_{n+1}$
and since $y\in P_n$ we have $f(y)\in P_{n+1}$.
Thus, $f(x),f(y)\in P_{n+1}$
and $\he (f(x)) =\he (f(y))$; so  $f(x) = f(y)$,
which is impossible since $f$ is one-to-one.

(g) Again, we have to prove that $f^n[S_0 ] =S_{n}$, for all $n\in \o$.
For $n=0$ this is trivial.
Assuming that $f^n[S_0 ] =S_{n}$, by (\ref{EQ8160'}) we have
$f^{n+1}[S_0 ]
=(f\circ f^{n})[S_0 ]
=f[f^{n}[S_0 ]]
=f[S_{n}]
=S_{n+1}$.

(h) Let $P_n=P_0$, where $n\in \N$.
Then, clearly, $S_n =S_0$ and we show that $S_0 =N_0\up$.
If $s\in S_0$, then $\he (s)\geq  \a  $
and, hence there is $u\in (\cdot ,s]\cap L_{\a } $.
For $t\in P_0$ we have $u,t \leq s$,
which implies $u \not\parallel t$
and, since $\he (t) < \a  $, we have $t<u$.
Thus $P_0 <u$ and, hence, $(\cdot ,u)=P_0 =(\cdot ,a_0)$,
that is $u\in N_0$ and $s\in N_0\up$.
Conversely, if $s\geq u\in N_0$, then $(\cdot ,u)=P_0$
and for $t\in P_0$ we have $t<u\leq s$ and, thus, $s\in S_0$.

By (g) we have $g:= f^n|S_0 \in \Cond (S_0, S_n)=\Cond (S_0) =\Cond (N _0 \up)$,
by (d) and (a) we have $b_0\in N_0$ and $\he (b_0) = \a   <\he (f(b_0))$
and, by Fact \ref{T8124}(a), $\he (f^n(b_0))\geq \he (f(b_0))$;
so $\he (g(b_0))> \he (b_0)$ and the node $N_0$ is critical indeed.

(i) By (a) and (d) we have $a_0, b_0 \in S_0$
and by (g) for each $n\in \N$ we have $a_n =f^n (a_0)\in f^n [S_0]=S_n$
and, similarly, $b_n \in S_n$.
Thus $P_n <T_n$, for all $n\geq 0$.
By (f) we have $f[S_{-1}]= S_0$ and, hence, $a_{-1}=f^{-1}(a_0)\in S_{-1}$
and, similarly, $b_{-1}\in S_{-1}$.
In the same way we show that $a_n,b_n\in S_{n}$, for all $n<0$.
So, $P_n <T_n$, for all $n< 0$, and (a3) is true.

Suppose that $P_m =P_n$, that is, $f^m [P_0]=f^n [P_0]$, for some $m<n\in \Z$.
Then $n-m\in \N$ and, by (b), $f^{n-m}[P_0]=P_0$,
that is $P_{n-m}=P_0$, which contradicts our assumption.
So, for different $m,n\in \Z$ we have $P_m\neq P_n$,
by (e), $P_n$, $n\in \Z$, are maximal chains of height $ \a  $ in the tree $\BT | \a  $.
Since $P_n <T_n$, for all $n\in \Z$, in the tree $\BT _f$ we have $\he (a_n),\he (b_n)\geq  \a  $.
So, the sets $P_n$, $n\in \Z$, are different maximal chains of height $ \a  $
in the tree $\BT _f| \a  = \bigcup _{n\in \Z} P_n$ and (a2) is true.

By (a) we have $a_1 <b_1$.
So, for each $n\in \N$, since $f^{n-1}$ is a homomorphism of $\BT$,
we have $a_n =f^{n-1}(a_1)< f^{n-1}(b_1)= b_n$.
So $\BT _n \cong 2$, for all $n>0$.
By (a) again, $a_0 \parallel b_0$, that is $f(a_{-1}) \parallel f(b_{-1})$,
which implies that $a_{-1} \parallel b_{-1}$
and, similarly, $a_{n} \parallel b_{n}$, for all $n<0$.
Thus $\BT _n$ is an antichain of size 2, for all $n\leq 0$, and (a4) is true.
This also implies that $\he (\BT _f)=  \a   +2$ and (a1) is true as well.

Since $f[T_f| \a  ]=f[\bigcup _{n\in \Z} P_n]=\bigcup _{n\in \Z} f[P_n]=\bigcup _{n\in \Z} P_{n+1}=T_f| \a   \subset T| \a  $,
from (b) it follows that $g:=f| (T_f| \a  )\in \Aut (\BT_f| \a  )$.
In addition we have $g[P_n]=f[f^n[P_0]]=f^{n+1}[P_0]=P_{n+1}$, for all $n\in \Z$,
and (a5) is true.

Since $a_{n+1}=f^{n+1}(a_0)=f(a_n)$ and, similarly, $b_{n+1}=f(b_n)$,
regarding (\ref{EQ000}) we have
$F_g := g\cup \{ \la a_n,a_{n+1}\ra :n\in \Z \}\cup \{ \la b_n,b_{n+1}\ra:n\in \Z \} =   f| (T_f| \a  )\cup \{ \la a_n,f(a_n) \ra :n\in \Z \}\cup \{ \la b_n,f(b_n)\ra:n\in \Z \}$
and, hence, $F_g =f| T_f :T_f \rightarrow T_f$ and the condensation $f$ extends $F_g$ to $T$.

(j)
First we prove that $S_n \cap L_\a =\{ s\in T : (\cdot ,s)=P_n \}$.
If $s\in S_n \cap L_\a$,
then $P_n <s$ and, hence, $P_n \subset (\cdot ,s)$.
Assuming that there is $t\in (\cdot ,s)\setminus P_n$
we would have $t\in T|\a$ and $P_n \cup \{ t\}$ would be a chain in $\BT |\a$,
which is false by (e).
Thus  $(\cdot ,s)=P_n$.
Conversely, if $(\cdot ,s)=P_n$,
then $s\in S_n$
and, by (e), $\he (s)=\otp (P_n)=\a$,
thus $s\in L_\a$.
The equality  is proved.

Clearly, $N_n =\{ s\in T : (\cdot ,s)=P_n \}\subset N$, for some $N\in \CN (\BT)$.
If $t\in N$ and $s\in N_n$,
then $(\cdot ,t)=(\cdot ,s)=P_n$
and, hence, $t\in N_n$. So, $N_n=N\in \CN (\BT)$.

If $t\in N_{n+1}$ then $t\in S_{n+1}$, $\he (t)=\a$
and, since by (f) $f |S_n\in \Cond (S_n ,S_{n+1})$, there is $s\in S_n$ such that $t=f(s)$.
Since $f\in \Cond (\BT )$ by Fact \ref{T8124}(a) we have $\a \leq \he (s) \leq \he (t)$,
which gives $s\in L_\a$ and, hence, $s\in N_n$.
Thus $t\in f[N_n]$ and the inclusion $N_{n+1}\subset f[N_n]$ is proved.
Since $f$ is a bijection we have  $|N_{n+1}|\leq |N_n|$.

(k) Let $N_0=\{ t_i:i<k\}$, where $k\in \N$, $t_0=b_0$ and $t_1=a_0$.
By (a) we have $f(b_0)\not\in L_\a $
and by (j) $N_1 \subset f[N_0]\cap L_\a \subset \{ f(t_i):0<i<k\}$,
which gives $|N_1|\leq k-1$.
The rest follows from (j) and (f).
\hfill $\Box$
\begin{rem}\label{R8104}\rm
If $\a$ from Proposition \ref{T8108}(a) is a successor ordinal, then there is $z_0\in L_{\a -1}$ such that $a_0,b_0 \in \is (z_0)$.
Defining $z_n =f^n(z_0)$, for $n\in \Z$, by (b) we have $z_n\in L_{\a -1}$
and by (e) $P_n =(\cdot , z_n]$ and $S_n =(z_n ,\cdot )=\is (z_n)\up$, for all $n\in \Z$.
\end{rem}
\begin{ex}\label{EX8104}\rm
{\it For each ordinal $\a>0$ the tree $\D_\a =\bcd _{n\in \Z}(\a + \BT _n)$ (see Example \ref{EX8101})
is an archetypical tree satisfying $\a _{\D _\a}=\a$. If $\a$ is a limit ordinal, then the restriction $\D _\a |\a _{\D _\a}$ is not reversible.}

Recall that $\D _\a=\bcd _{n\in \Z}(\P_n +\BT _n)$,
where $\P_n$, $n\in \Z$, are pairwise disjoint copies of $\a$,
$T _n =\{ a_n, b_n\}$, $a_n \parallel b_n$, for $n\leq 0$, and $a_n < b_n$, for $n> 0$.
For the mapping $g$ satisfying (a5) and $F_g$ defined by (\ref{EQ000}) we have $\a _{F_g} =\a$; thus $\a _{\D _\a}\leq\a$.

First suppose that $\a _{\D _\a}=0$.
Let $f \in \Cond (\D _\a )\setminus \Aut (\D _\a )$, where $\a _f=0$,
and let $t_0 \in L_0$ be such that $\he (f(t_0))>\he (t_0) (=0)$.
Then there are $m,n\in \Z$ such that $t_0\in P_m$ and $f(t_0)\in P_n \cup T_n$.
For each $t\in P_m \cup T_m$ we have $t_0 \leq t$ and, hence,  $f(t_0) \leq f(t)$,
which gives $f[P_m \cup T_m]\subset [f(t_0), \cdot) \subset P_n \cup T_n$.
For $t\in T_m$ we have $\he (t)\geq \a$ and, hence, $\he (f(t))\geq \a$;
thus $f[T_m]\subset T_n$ and, since $|T_m|=|T_n|=2$, $f[T_m]=T_n$.
Clearly we have $s:=\min P_n <f(t_0)$
and, since $f$ is onto, there is $s_0 \in T\setminus (P_m \cup T_m)$ such that $f(s_0)=s$.
Since $\he (s_0)\leq \he (f(s_0))=0$ we have $s_0\in L_0$
and there is $k\in Z \setminus \{ m\}$ such that $s_0=\min P_k$.
Now as above we prove that $f[T_k]=T_n$,
which is impossible because $f[T_m]=T_n$ and $f$ is one-to-one.

Thus $\a _{\D_\a} >0$.
Let $f \in \Cond (\D _\a )\setminus \Aut (\D _\a )$, where $\a _f=\a _{\D_\a}$.
By (a) and (d) of Proposition \ref{T8117} there are $a,b\in L_{\a _{\D_\a}}$
belonging to the same node $N$ of $\D_\a$.
Since $a\parallel b$, this is possible only if $\{ a, b\}=\{ a_n, b_n\}$, for some $n\leq 0$.
Since $a\in L_{\a _{\D_\a}}$ and $a_n\in L_\a$ we have $\a _{\D_\a}=\a$.

If $\a$ is a limit ordinal, then $\D _\a |\a _{\D _\a}=\bcd _{n\in \Z}\a$ is not reversible (see Fact \ref{T000}).
\end{ex}
\section{Trees with finite nodes. Regular $\l$-ary trees}\label{S5}
The following theorem, based on Proposition \ref{T8108}, provides detection of some classes of reversible trees with finite nodes and several
constructions of such trees.
\begin{te}\label{T8118}
If $\BT$ is a tree with finite nodes, then each of the following conditions implies that $\BT $ is reversible

(i) $\he (\BT)\leq \o +1$,

(ii) For each ordinal $\a \in [\o ,\he (\BT ))$
the sequence $\la |N|: \CN (\BT ) \ni N\subset L_\a\ra $ is almost constant or the sequence $\la\la |N|,|N\up|\ra : \CN (\BT ) \ni N\subset L_\a\ra $ is finite-to-one.
\end{te}
\dok
We prove the contrapositive.
Let $\BT$ be a non-reversible tree with finite nodes.
By Proposition \ref{T8117}(a) $\BT$ has no critical nodes
and by Theorem \ref{T8198} it contains an archetypical tree.
Taking $f\in \Cond (\BT)\setminus \Aut (\BT)$ and using Proposition \ref{T8108}
we prove that $\neg $(i) and $\neg $(ii) hold. Let $\a  :=\min \{ \he (t): t\in T \land \he (t)< \he (f(t))\}$.

Since $L_0 \in \CN (\BT)$ we have $|L_0|<\o$
and an easy induction shows that $|L_n|<\o$, for all $n\in \o$.
Since $P_n$, $n\in \Z$, are different branches in $\BT |\a$, we have $\a \geq \o$.
In addition, $\he (f(a))>\a$, thus $\he (f(a))\geq \o +1$ and, hence, $\he (\BT )> \o +1$.
So $\neg $(i) is true.

Since $|N_0|<\o$,
by (j) and (k) of Proposition \ref{T8108} we have
$\dots\geq|N_{-2}|\geq |N_{-1}|\geq|N_0|>|N_1|\geq |N_2|\geq \dots$,
thus the sequence $\la |N|: \CN (\BT ) \ni N\subset L_\a\ra $ is not almost constant.
In addition, there is $n_0\in \N$ such that  $|N_n|=|N_{n_0}|$ and $|N_n \up|=|N_{n_0} \up|$, for all $n\geq n_0$,
and  $\la\la |N|,|N\up|\ra : \CN (\BT ) \ni N\subset L_\a\ra $ is not a finite-to-one sequence.
So, $\neg $(ii) is true as well.
\hfill $\Box$
\begin{ex}\rm
{\it Constructions of reversible trees.}
First, let $\N =C\du F$
and, for $k\in \N$, let $\la n^k_i:i \in \o\ra \in \N ^\o$ be
an almost constant sequence, if $k\in C$;
or a finite-to-one sequence, if $k\in F$.
We construct a countable regular reversible tree $\BT ={}^{<\o}2 \cup \bigcup _{k\in \o }L_{\o +k}$
of height $\o +\o$,
defining the levels $L_{\o +k}$, for $k\in \o$, by recursion.
First, let $x_i$, $i\in \o$, be different elements of ${}^{\o}2$ and $L_\o =\{ x_i :i\in \o \}$.
If $L_{\o +k}$ is defined and $L_{\o +k}=\{ t_i :i\in \o \}$ is an enumeration,
then to each $t_i$ we add $n_i$-many immediate successors and define $L_{\o +k+1} :=\bigcup _{i\in \o}\is (t_i)$.
Thus,  for each $k\in \o$ we have $\la |N|: \CN (\BT ) \ni N\subset L_{\o +k+1}\ra=\la n^k_i:i \in \o\ra $
(up to an re-enumeration of the sequence)
and condition (ii) of Theorem \ref{T8118} is satisfied.
In this way we obtain $\c$-many non-isomorphic countable reversible trees.
Clearly there are several variations of the above construction
(e.g.\ producing trees of bigger height, size or with maximal elements) with the restriction that
the sequence $\la |N|: \CN (\BT ) \ni N\subset L_\a\ra $ must be almost constant, if the level $L_\a$ is uncountable.

Another simple construction is to
take countably many different elements of ${}^{\o}2$ indexed by $\N ^2$,
define $L_\o =\{ x_{m,n} :\la m,n\ra\in \N^2 \}$,
for each $\la m,n\ra\in \N^2$ add a copy $T_{m,n}$ of the tree $\bcd _m n$ above $x_{m,n}$
and define $\BT ={}^{<\o}2 \cup L_\o \cup \bigcup _{\la m,n\ra\in \N^2}T_{m,n}$.
Then the sequence $\la |N|: \CN (\BT ) \ni N\subset L_{\o +k}\ra $ is constant at 1, for each $k\in \o \setminus \{1\}$,
and it is neither almost constant nor finite-to-one for $k=1$.
But  $\la\la |N|,|N\up|\ra : \CN (\BT ) \ni N\subset L_{\o+1}\ra =\la\la m,mn\ra : \la m,n\ra\in \N^2\ra$
is an one-to-one sequence and the tree $\BT$ is reversible by Theorem \ref{T8118}.

Taking $\bcd _m \o _n$ instead $\bcd _m n$ we obtain a reversible tree again
and, in addition, the sequence $\la |N\up|: \CN (\BT ) \ni N\subset L_{\o +1}\ra= \la \o _n : \la m,n\ra\in \N^2\ra$
is neither almost constant nor finite-to-one.
\end{ex}
\begin{te}\label{T8119}
(a) For each $n\in \N$, each regular $n$-ary tree is reversible.
In particular, for each ordinal $\a >0$ the complete $n$-ary tree of height $\a$, ${}^{<\a }n$, is reversible.

(b) There is a reversible $n$-ary Aronszajn tree, for each integer $n\geq 2$.

(c) If there exists a Suslin or a Kurepa tree, then  there is a reversible $n$-ary one, for each integer $n\geq 2$.
\end{te}
\dok
If $\BT$ is a regular $n$-ary tree, $\a \in [\o ,\he (\BT))$ and $\CN (\BT)\ni N\subset L_\a$,
then $|N|=n$, if $\a$ is a successor ordinal, and  $|N|=1$, otherwise.
So, by Theorem \ref{T8118} the tree $\BT$ is reversible and (a) is true.

Statements (b) and (c) are true because by Proposition 2.8 of \cite{Tod} if $\BT$ is an Aronszajn or Suslin or Kurepa tree,
then for each integer $n\geq 2$ the tree $\BT$ contains a regular $n$-ary subtree $\BT'$ of the same type, which is  reversible by (a).
\kdok
The following Proposition  \ref{T8121} describes one type of a critical node and will be used in a characterization of reversible regular $\l$-ary trees.
First, by Theorem 3.2 of \cite{KuMo4} we have
\begin{fac}[\cite{KuMo4}]\label{TA019}
Let $\X _i$, $i\in I$, be pairwise disjoint and connected  binary structures.
Then the structure $\bcd _{i\in I} \X _i$ is not reversible  iff
there are a surjection  $f:I\rightarrow I$ and monomorphisms
$g_i :\X _i \rightarrow \X _{f(i)}$, for $i\in I$, such that $\{g_i[X_i]: i\in f^{-1}[\{ j\}] \}$ is a partition of the set $X_j$,
for each $j\in I$, and the function  $f$ is not one-to-one or  $g_i \not\in \Iso (\X _i,\X _{f(i)})$, for some $i\in I$.
\end{fac}
\begin{prop}\label{T8121}
If $\BT$ is a tree, $N\in \CN (\BT )$, $s_i$, $i\in \o$, are different elements of $N$, $T'\subset [s_0,\cdot)$ and there are condensations
$g_0 : [s_0 ,\cdot ) \rightarrow [s_0 ,\cdot )\setminus T',$
$g_1 : [s_1 ,\cdot ) \rightarrow T'$ and
$g_i : [s_i ,\cdot ) \rightarrow [s_{i-1} ,\cdot )$, for $i\geq 2$,
then the tree $\BT$ is not reversible. In addition, $N$ is a critical node of $\BT$.
\end{prop}
\dok
Since $S=\{ s_i : i\in \o\}\subset N$,
by Fact \ref{T8111}(a) it is sufficient to prove that the tree $S\up =\bcd _{i\in \o } [s_i ,\cdot )$ satisfies the assumptions of Fact \ref{TA019}.
First, the trees $\X _i:=[s_i ,\cdot )$, $i\in \o$, are pairwise disjoint, rooted and, hence, connected binary structures.
If $f:\o \rightarrow \o$, where $f(0)=0$ and $f(i)=i-1$, for $i\geq 1$,
then $f\in \Sur (\o )\setminus \Sym (\o )$ and, by our assumption,  $g_i \in \Mono (X_i , X_{f(i)})$, for each $i\in \o$.
In addition, $\{g_i[X_i]: i\in f^{-1}[\{ 0\}] \}=\{ [s_0 ,\cdot )\setminus T' , T'\}$ is a partition of the set $X_0$
and, since $f^{-1}[\{ j\}]=\{ j+1\}$ and  $g_{j+1} [[s_{j+1} ,\cdot )] = [s_j ,\cdot )$, for $j>0$,
the assumptions of Fact \ref{TA019} are satisfied and the tree $S\up$ is not reversible.
It is easy to check that $g=\bigcup _{i\in \o}g_i \cup \id _{N\up \setminus S\up}\in \Cond (N\up)$. If $s_0\in T'$, then $\he (g(s_0))>\he (s_0)$;
otherwise we have $\he (g(s_1))>\he (s_1)$ and $N$ is a critical node of $\BT$ indeed.
\hfill $\Box$
\begin{ex}\label{EX8103}\rm
Homogeneous $\l$-ary trees are non-reversible, for $\l \geq \o$.
Defining a tree $\BT$ to be {\it homogeneous} iff $\BT \cong [t, \cdot)$, for all $t\in T$,
we easily check that a homogeneous tree $\BT$ must be rooted, $\l$-ary, for some cardinal $\l >0$ and of height $\o ^\d$, for some ordinal $\d$.
(We recall that the ordinals of the form $\o ^\d$ are called {\it indecomposable}.)
If, in addition, $|\l |\geq \o$ and, hence, $\d >0$, then  taking $\{ s_i : i\in \o\}\subset L_1=:N$ and $t\in \is (s_0)$,
we easily find condensations (in fact, isomorphisms) $g_i$, $i\in \o$, satisfying the assumptions of Proposition \ref{T8121}.
\end{ex}
\begin{te}\label{T8127}
The tree  ${}^{<\a }\l$ is reversible  iff $\,\min \{\a ,\l \}<\o$.
\end{te}
\dok
The implication ``$\Leftarrow$" follows from  Theorem  \ref{T8197}(a) and Theorem \ref{T8119}(a).

For the converse we suppose that $\a ,\l \geq \o$ and, using Proposition \ref{T8121}, prove that $\BT :={}^{<\a }\l$ is a non-reversible tree.
First, since $\a \geq \o$, for each $s=\la \xi\ra\in L_1$ we have $\BT \cong [s ,\cdot )$;
one isomorphism $F_s:\BT \rightarrow [s ,\cdot )$ is defined by $F_s(t)=\la\xi\ra^{\smallfrown}t =\la \xi ,t (0),t (1),\dots \ra$, for all $t \in {}^{<\a }\l$.
More precisely, if $\eta :\N \rightarrow \o$, where $\eta (n)=n-1$, for each $n\in \N$,
then $F_s(t)=\{\la 0,\xi\ra\} \cup ((t \upharpoonright \o) \circ \eta )\cup (t \upharpoonright (\a \setminus \o))$.

Second, since $\l\geq \o$, for each $s\in L_1$ we have $\BT \cong T\setminus [s ,\cdot )$, because both trees are rooted unions of $\l$-many copies of $\BT$

Let $s_i$, for $i\in\o$, be different elements of the node $L_1$ of $\BT$ and, say, $s_0=\la \xi\ra$.
Then $t:=F_{s_0}(s_0)=\la \xi ,\xi\ra\in \Lev_1 [s_0 ,\cdot )$
and, as above, there is an isomorphism $g_0:[s_0, \cdot )\rightarrow [s_0, \cdot )\setminus [t ,\cdot )$
Since $[t ,\cdot )=F^2[T]\cong \BT\cong [s_1, \cdot )$
there is an isomorphism $g_1 :[s_1, \cdot )\rightarrow [t ,\cdot )$
and, clearly, there are isomorphisms $\;g_i : [s_i ,\cdot ) \rightarrow [s_{i-1} ,\cdot )$, for $i\geq 2$.
By Proposition \ref{T8121} the tree  ${}^{<\a }\l$ is not reversible.
\hfill $\Box$
\section{Unions of complete $\l$-ary trees}\label{S6}
By Fact \ref{T8111}(a) in the analysis of the trees having infinite nodes
the results concerning reversibility of disjoint unions of rooted trees can help.
One step towards the understanding the general situation
is to find out what is going on with the trees of the form $\bcd _{i\in I}{}^{< \a _i}\l _i$,
where $\l _i>0$ are cardinals and $\a  _i>0$ are ordinals,
and Fact \ref{T000} shows that even when $\l _i=1$, for all $i\in I$,
a combinatorial characterization of reversibility is non-trivial.

Generally speaking, in several classes of trees
extreme properties provide reversibility;
for example, both ``similarity" (almost constant)
and ``diversity" (finite-to one) in Theorem \ref{T8118},
or both ``low and wide" and "tall and narrow" in Theorem \ref{T8127}.
In the sequel we show that similar phenomena appear in our case.

By Fact \ref{TB051} and Theorem \ref{T8127}
a necessary condition for the reversibility of the tree $\bcd _{i\in I}{}^{< \a _i}\l _i$
is that $\a _i <\o$ or $\l _i <\o$, for all $i\in I$,
and by Fact \ref{T001}
it is sufficient that the sequence $\la \la |{}^{< \a _i}\l _i|, \a _i\ra :i\in I\ra$ is finite-to-one
(for example, the tree $\bigcup _{\la m,n\ra \in [2, \o)\times \o}{}^{<m}\o _n$ is reversible).
First we regard the simple situation, when $\BT =\bcd _\mu {}^{<\a }\l $,
where $\mu >0$ is a cardinal.
For $\l =1$ we have $\BT =\bcd _\mu \a$
and by Fact \ref{T000} $\BT$ is reversible
iff $\mu <\o$ or $\a$ is a successor ordinal.
So, in the sequel we consider the case when $\l \geq 2$ and first prove an auxiliary statement.
\begin{lem}\label{T8128}
If $\BT $ is a tree, $\d +1 <\a =\he (\BT)$, $G\in \Cond (\BT |(\d +1) )$, $G[L_\d]=L_\d$ and
$f_t\in \Cond ([t, \cdot ), [G(t), \cdot ))$, for $t\in L_\d$, then $F:= G\cup \bigcup _{t\in L_\d}f_t \in \Cond (\BT )$.
\end{lem}
\dok
For $t\in L_\d$ we have $f_t\in \Cond ([t, \cdot ), [G(t), \cdot ))$ and, hence, $f_t(t)=G(t)$.
Since $G|L_\d : L_\d \rightarrow L_\d$ is a bijection, $\{[t, \cdot ):t\in L_\d \}=\{[G(t), \cdot ):t\in L_\d \}$ is a family of pairwise disjoint trees
and $f:=\bigcup _{t\in L_\d}f_t \in  \Aut(\bigcup _{t\in L_\d}[t, \cdot ))$.
In addition we have $G|L_\d =f|L_\d$ and $F$ is a well defined bijection.
For a proof that $F$ is a homomorphism, for $a,b\in T$, where $a<b$, we show that $F(a)<F(b)$.
Since $F\upharpoonright (\BT |(\d +1))= G \upharpoonright(\BT |(\d +1))$  and $F\upharpoonright \bigcup _{t\in L_\d}[t, \cdot )=f \upharpoonright \bigcup _{t\in L_\d}[t, \cdot )$
the only non-trivial case is when $a\in \BT |\d$ and $b\in \bigcup _{t\in L_\d}(t, \cdot )$.
Then we have $a,t <b$ and, hence, $a<t<b$,
which gives $F(a)=G(a)<G(t)=f(t)<f(b)=F(b)$ and we are done.
\hfill $\Box$
\begin{te}\label{T8126}
Let $\l \geq 2$. The tree $\BT :=\bcd _\mu {}^{<\a }\l$ is reversible  $\Leftrightarrow\min\{\a ,\l\mu\} <\o$.
\end{te}
\dok
($\Leftarrow$)
If $\a =n\in \N$,
then $\he (\BT )=n$
and all branches of  $\BT$ are of height $n$;
so, by  Theorem \ref{T8197}(a), $\BT$ is reversible.
If $\l=n\in \N$ and $\mu =m\in \N$,
then by Theorem \ref{T8119}(a) the connectivity components of $\BT$, $\BT_i \cong {}^{<\a }n$, $i<m$, are reversible
and $\BT$ is reversible by Fact \ref{TB051}(b).

($\Rightarrow$) Assuming that $\a \geq \o$ and $\l\mu \geq \o$, we show that $\BT$ is not reversible.
If $\l \geq \o$, then by Theorem \ref{T8127} the tree ${}^{<\a }\l$ is not reversible and $\BT$ is not reversible by Fact \ref{TB051}(a).
Thus it remains to be proved that $\BT$ is not reversible when $2\leq \l =k\in \N$ and $\mu \geq \o$,
moreover, by Fact \ref{TB051}(a), it is sufficient to prove the statement when $\mu =\o$, that is
$$\textstyle
\BT =\bcd _\o {}^{<\a }k.
$$

{\it Case 1}: $\a =\o$. Then $\BT =\bcd _\o {}^{<\o }k$ is a bad tree of height $\o$ and it is not reversible by Theorem \ref{T8110}(a).

{\it Case 2}: $\a =\o +1$. Then $\BT =\bcd _\o {}^{\leq \o }k=\bcd _{i<\o}\BT _i$.
First we consider the tree $\BT _0$; say $\BT_0={}^{\leq\o }k$.
If $\f \in {}^{\leq\o }k$, then
either $\ran \f \subset \{ 0\}$, that is $\f =0_\a$, for some $\a \leq \o$ (so, $0_0=\emptyset$ and $0_\o$ is the $\o$-sequence of zeroes)
or $\ran \f \not\subset \{ 0\}$, which means that $\ran \f\cap \{1,\dots ,k-1\}\neq \emptyset$.
Let $S_0$ be the subtree of the tree ${}^{\leq\o }k$ defined by $S_0 := \{ \emptyset \} \cup [0 ,\cdot )$.
Then $S_0=\{ 0_\a :\a \leq \o\}\cup \{ 0^\smallfrown \f:\f \in {}^{\leq\o }k \land \ran \f \not\subset \{ 0\} \}$;
clearly, $0^\smallfrown \f$ is the sequence $\la 0, \f (0), \f (1), \dots\ra=\{ \la 0,0\ra\}\cup (\f \circ \eta)$,
where $\eta :\N \rightarrow \o$ is the bijection given by $\eta (n)=n-1$.

Further, we prove that the mapping $g_0:T_0\rightarrow S_0$ is a condensation, where
\begin{equation}\label{EQ8166}
g_0(\f )= \left\{
                   \begin{array}{ll}
                      \f ,              & \mbox { if } \; \ran \f \subset \{ 0\}, \\
                      0^\smallfrown \f, & \mbox { if } \; \ran \f \not\subset \{ 0\} .
                   \end{array}
            \right.
\end{equation}
First we show that $g_0 [{}^{\leq\o }k]= S_0$.
Let $\f \in  {}^{\leq\o }k$. If $\f =0_\a$, for some $\a\leq \o$, then $g_0(\f) =\f =0_\a \in S_0$;
otherwise we have $g_0(\f) =0^\smallfrown \f\in [0 ,\cdot ) \subset S_0$ and, thus,  $g_0[{}^{\leq\o }k]\subset S_0$.
Let $\p \in S_0$.
If $\p =0_\a$, for some $\a\leq \o$, then $\p =g_0(\p )\in g_0[{}^{\leq\o }k]$.
Otherwise, $\p =0^\smallfrown \f$, where $\f \in {}^{\leq \o }k$ and $\ran \f  \not\subset \{ 0\}$, so by (\ref{EQ8166}) we have $\p =g_0(\f )\in g_0[{}^{\leq\o }k]$ again.
Thus, $g_0[{}^{\leq\o }k]\supset S_0$ and  $g_0$ is a surjection.

By (\ref{EQ8166}), for $\f \in {}^{\leq\o }k$ we have: $\ran \f\subset \{ 0\}$ iff $\ran g_0(\f )\subset \{ 0\}$.
Let $g_0(\f _1 )=g_0 (\f _2)$.
If $\ran g_0 (\f _1 )\subset \{ 0\}$, then $\ran \f _1 , \ran \f _2 \subset \{ 0\}$ and, by (\ref{EQ8166}), $\f _1 =g_0(\f _1)=g_0(\f _2 )=\f _2$.
If $\ran g_0(\f _1 )\not\subset \{ 0\}$, then $\ran \f _1 , \ran \f _2 \not\subset \{ 0\}$
and, by (\ref{EQ8166}) again, $0^\smallfrown\f _1 =g_0 (\f _1)=g_0 (\f _2 )=0^\smallfrown\f _2$,
which implies $\f _1 =\f _2$, so $g_0$ is an injection.

Let $\f _1, \f _2 \in {}^{\leq\o }k$ and $\f _1\varsubsetneq \f _2$; then $|\f _1|<\o$.
If $\ran \f _2 \subset \{ 0\}$, then $\ran \f _1 \subset \{ 0\}$,
which by (\ref{EQ8166}) gives $g_0(\f _1) =\f _1\varsubsetneq \f _2 =g_0(\f _2)$.
If $\ran \f _1 \not\subset \{ 0\}$, then $\ran \f _2 \not\subset \{ 0\}$,
which implies that $g_0(\f _1) =0^\smallfrown\f _1\varsubsetneq 0^\smallfrown\f _2 =g_0(\f _2)$.
Finally, if $\ran \f _1 \subset \{ 0\}$ and $\ran \f _2 \not\subset \{ 0\}$,
then, since  $|\f _1|<\o$, there is $n\in \o$ such that $\f _1 =0_n$
and, hence, $\f _2 = 0_n ^\smallfrown \p$, where $\ran \p\not\subset \{ 0\}$.
Now,  by (\ref{EQ8166}) we have $g_0(\f _1) =0_n \varsubsetneq 0 ^\smallfrown 0_n ^\smallfrown \p = g_0(\f _2)$
and $g_0$ is a homomorphism. Thus $g_0$ is a condensation indeed.

Let $\{ S_i:i<k\}$ be a partition of the tree $\BT _0 ={}^{\leq \o }k$, where the set $S_0$ is defined above and $S_i := [i ,\cdot )_{\BT _0}$, for $i\in [1,k)$.
Clearly for $i\in \{1,\dots ,k-1\}$ we have $S_i \cong {}^{\leq \o }k \cong \BT _i$;
let $g _i :\BT _i\rightarrow S_i$, for $i\in \{1,\dots ,k-1\}$ be isomorphisms.
Thus we have monomorphisms $g_i : \BT _i \rightarrow \BT _0$, for $i<k$,
such that $\{ g_i [T_i]: i<k\}$ is a partition of the set $T_0$.
Taking isomorphisms $g_i:\BT _i\rightarrow \BT _{i-k+1}$, for $k\leq i<\o$,
and defining $f:\o \rightarrow \o$ by $f(i)=0$, for $i<k$, and $f(i)=i-k+1$, for $k\leq i<\o$,
we have $f\in \Sur (\o )\setminus\Sym (\o )$.
Now the tree $\BT$ is not reversible by Theorem \ref{TA019} and $G=\bigcup _{i<\o}g_i\in \Cond (\BT )\setminus \Aut (\BT )$.

{\it Case 3}: $\a >\o +1$.
If $G$ is the condensation of the restriction $\BT |(\o +1)\cong \bcd _\o {}^{\leq \o }k$ from Case 2,
then by the construction we have $G\not\in \Aut (\BT |(\o +1))$ and
$G [\Lev _\o (\BT |(\o +1))]=\Lev _\o (\BT |(\o +1))$.
In addition, for each $t\in \Lev _\o (\BT |(\o +1))$ we have $[t,\cdot )\cong {}^{<\b }k$, where $\b \cong\a \setminus(\o +1)$,
and taking isomorphisms $f_t : [t,\cdot )\rightarrow [G(t),\cdot )$, for $t\in \Lev _\o (\BT |(\o +1))$, and defining
$F=G\cup \bigcup _{t\in \Lev _\o (\BT |(\o +1))}f_t$, by Lemma \ref{T8128} we have $F\in \Cond (\BT )\setminus \Aut (\BT )$.
So the tree $\BT =\bcd _\o {}^{<\a }k$ is not reversible indeed.
\kdok
If $\l \geq \o$,
then by Theorem \ref{T8126} the tree $\bcd _\mu {}^{<\a}\l$ is reversible iff its height $\a$ is finite.
In the sequel we regard the unions of $\l$-ary trees of different height of the form $\bcd _{i\in I}{}^{<n_i}\l$, where $\l \geq \o$, and characterize the reversible ones in that class.
\begin{lem}\label{T100}
If $0<m<n<\o$ and $\l$ is an infinite cardinal,
then there are embeddings $f: {}^{<m}\l \rightarrow {}^{<n}\l$ and $g: {}^{<n}\l \rightarrow {}^{<n}\l$
such that $\{ f[{}^{<m}\l], g [{}^{<n}\l]\}$ is a partition of the tree ${}^{<n}\l$.
\end{lem}
\dok
For $k\in \N$ let $0^k$ denote the sequence $\la 0,\dots,0\ra$ with $k$ zeroes.
It is easy to check that the mappings $f: {}^{<m}\l \rightarrow [0^{n-m}, \cdot)_{{}^{<n}\l}$,
where $f(\f )=0^{n-m}{}^{\smallfrown}\f$, for all $\f\in {}^{<m}\l$, and $g: {}^{<n}\l \rightarrow {}^{<n}\l \setminus [0^{n-m}, \cdot)_{{}^{<n}\l}$, defined by
$$
g(\f ) = \left\{
                     \begin{array}{ll}
                      \f , & \mbox { if }  \f \in {}^{<n}\l \setminus \bigcup _{k\in \o}[0^{n-m-1}{}^{\smallfrown}k, \cdot), \\[2mm]
                      0^{n-m-1}{}^{\smallfrown}(k+1)^{\smallfrown}\p,      & \mbox { if }  \f =0^{n-m-1}{}^{\smallfrown}k^{\smallfrown}\p, k\in \o \mbox{ and }\p\in {}^{<m}\l,
                     \end{array}
                   \right.
$$
are isomorphisms.
\hfill $\Box$
\begin{te}\label{T101}
If $\l\geq \o$, $\la n_i :i\in I\ra \in \N ^I$ and $I_k := \{ i\in I :n_i=k\}$, for $k\in \N$,
then the tree $\BT=\bcd _{i\in I}{}^{<n_i}\l$ is reversible iff the sequence $\la n_i :i\in I\ra $ is

(i)  either finite-to-one,

(ii) or almost constant at some $k\in \N$ and bounded by $k$; that is,
\begin{equation}\label{EQ100}
\exists k\in \N \;\;\Big(|I\setminus I_k|<\o \land \forall i\in I\setminus I_k \;\;n_i <k \Big).
\end{equation}
\end{te}
\dok
First, by Theorem \ref{T8197}(a) the trees $\BT _i:={}^{<n_i}\l$, for $i\in I$, are reversible.

($\Leftarrow$) If (i) holds, then, since $\he (\BT _i)=n_i$, $\BT$ is reversible by Fact \ref{T001}.

Let (ii) hold and  $f\in \Cond(\BT )$.
For $i\in I$ let $r_i$ be the root of $\BT _i$.
If $i\in I\setminus I_k$, then $r_i =f (r_j)$, for some $j\in I$,
and assuming that $j\in I_k$,
there would be a $k$-sized chain $C$ in $\BT_j$ and $r_j \in C$;
so $f[C]$ would be a  $k$-sized chain in $\BT_i$, which is impossible.
Thus $j\in I\setminus I_k$ and we have proved that $\{ r_i : i\in I\setminus I_k\}\subset \{ f(r_i) : i\in I\setminus I_k\}$.
Since these sets are finite we have  $\{ r_i : i\in I\setminus I_k\}= \{ f(r_i) : i\in I\setminus I_k\}$.
So there is a permutation $\pi \in \Sym (I\setminus I_k)$ such that $f(r_i)=r_{\pi (i)}$
and, hence, $f[T_i]\subset T_{\pi (i)}$, for all $i\in I\setminus I_k$,
which gives $\bigcup _{i\in I\setminus I_k}f[T_i]\subset \bigcup _{i\in I\setminus I_k}T_i$.
Assuming that there is $t\in \bigcup _{i\in I\setminus I_k}T_i \setminus \bigcup _{i\in I\setminus I_k}f[T_i]$
there would be $i\in I\setminus I_k$, $j\in I_k$ and $s\in T_j$ such that $T_i\ni t=f(s)$
and, as above, there would  be a  $k$-sized chain in $\BT_i$, which is false.
Thus $f[\bigcup _{i\in I\setminus I_k}T_i]=\bigcup _{i\in I\setminus I_k}f[T_i]= \bigcup _{i\in I\setminus I_k}T_i$
and, hence,  $f[\bigcup _{i\in I_k}T_i]= \bigcup _{i\in  I_k}T_i$.
By Fact \ref{TB051}(b) and Theorem \ref{T8197}(a) the trees $\bigcup _{i\in I\setminus I_k}\BT _i$ and $\bigcup _{i\in I_k}\BT _i$ are reversible;
so, the corresponding restrictions of $f$ are automorphisms,  which implies that $f\in \Aut (\BT )$.

($\Rightarrow$) We prove the contrapositive: assuming $\neg\,$(i) and $\neg\,$(ii) we show that $\BT$ is not reversible.
By $\neg\,$(i) we have $K:=\{ k\in \N : |I_k|\geq \o\}\neq \emptyset$.

First we show that there are $k\in K$ and $i\in I$ such that $n_i >k$.
If $|K|>1$, this is evident; so, let $K=\{ k\}$; then $|I_k|\geq \o$.
If $|I\setminus I_k|<\o$, then by $\neg\,$(ii) there is $i\in I\setminus I_k$ such that $n_i>k$.
If $|I\setminus I_k|\geq\o$,
then, since $|K|=1$, the sequence $\la n_i :i\in I\setminus I_k\ra$ is finite-to-one
and, hence, there is $i\in I\setminus I_k$ such that $n_i>k$.

So, let $k\in K$, $i_0\in I$ and $n_{i_0} >k$.
Let us take different $i_j\in I_k$, for $j\in \N$.
Then $S:=\{ r_{i_j}:j\in \o\}\subset L_0 \in \CN (\BT)$
and by Fact \ref{T8111}(a) it is sufficient to prove
that the tree $S\up =\bigcup _{j\in \o}[r_{i_j}, \cdot )=\bigcup _{j\in \o}\BT_{i_j}$ is not reversible.
Since $[r_{i_0}, \cdot )\cong {}^{<n_{i_0}}\l$ and $[r_{i_j}, \cdot )\cong {}^{<k}\l$, for $j>0$,
by Lemma \ref{T100} there are embeddings $g_0:[r_{i_0}, \cdot )\rightarrow [r_{i_0}, \cdot )$
and $g_1:[r_{i_1}, \cdot )\rightarrow [r_{i_0}, \cdot )$
such that $\{ g_0[[r_{i_0}, \cdot )],g_1[[r_{i_1}, \cdot )]\}$  is a partition of $[r_{i_0}, \cdot )$.
Taking isomorphisms $g_j :[r_{i_j}, \cdot )\rightarrow [r_{i_{j-1}}, \cdot )$, for $j\geq 2$,
and using Proposition \ref{T8121} we conclude that the tree $S\up$ is not reversible.
\hfill $\Box$
\paragraph{Acknowledgement.}
This research was supported by the Science Fund of the Republic of Serbia,
Program IDEAS, Grant No.\ 7750027:
{\it Set-theoretic, model-theoretic and Ramsey-theoretic
phenomena in mathematical structures: similarity and diversity}--SMART.

\end{document}